\documentclass[12pt]{amsart}

\usepackage{amsmath}
\usepackage{amssymb}
\usepackage{latexsym}
\usepackage{amsfonts}

\usepackage[pagebackref,hypertexnames=false, colorlinks, citecolor=black, linkcolor=blue, urlcolor=red]{hyperref}

\input txdtools

\newcommand{\be}{\begin{equation}}
\newcommand{\ee}{\end{equation}}

\numberwithin{equation}{section}

\newtheorem{thm}{Theorem}[section]

\newtheorem{cor}[thm]{Corollary}

\newtheorem*{prop*}{Proposition}

\theoremstyle{remark}

\newtheorem*{rem*}{Remark}

\begin{document}

\title[Matrix Carleson Embedding and Sparse Operators]{A Study of the Matrix Carleson Embedding Theorem with Applications to Sparse Operators}
\date{\today}

\author[K. Bickel]{Kelly Bickel$^{\dagger}$}
\thanks{$\dagger$ Research supported in part by National Science Foundation grant DMS \# 1448846}
\address{Kelly Bickel, Department of Mathematics\\
Bucknell University\\
701 Moore Ave\\
Lewisburg, PA 17837}
\email{kelly.bickel@bucknell.edu}

\author[B. D. Wick]{Brett D. Wick$^{\ddagger}$}
\thanks{$\ddagger$ Research supported in part by National Science Foundation grant DMS \# 0955432}
\address{Brett D. Wick, School of Mathematics\\
Georgia Institute of Technology\\
686 Cherry Street\\
Atlanta, GA USA 30332-0160}
\email{wick@math.gatech.edu}

\begin{abstract}
In this paper, we study the dyadic Carleson Embedding Theorem in the matrix weighted  setting.  We provide two new proofs of this theorem, which highlight connections between the matrix Carleson Embedding Theorem and both maximal functions and $H^1$-BMO duality. Along the way, we establish boundedness results about new maximal functions associated to matrix $A_2$ weights and duality results concerning $H^1$ and BMO sequence spaces in the matrix setting. As an application, we then use this Carleson Embedding Theorem to show that if $S$ is a sparse operator, then the operator norm of $S$ on $L^2(W)$ satisfies:
\[ \|  S\|_{L^2(W) \rightarrow L^2(W)} \lesssim [W]_{A_2}^{\frac{3}{2}},\]
for every matrix $A_2$ weight $W$.
\end{abstract}

\maketitle

\section{Introduction}
The subject of this paper is the dyadic Carleson Embedding Theorem - an important tool for establishing the boundedness of paraproducts via testing conditions - and its applications and shortcomings in the matrix weighted setting.  This paper is particularly motivated by interest in the matrix $A_2$ conjecture. To set the scene, let us briefly review the scalar situation. 

\subsection{Scalar $A_2$ Conjecture.} Let $T$ be a Calder\'on-Zygmund operator and $w(x)$ a weight, i.e.~a locally integrable function on $\mathbb{R}$ that is positive almost everywhere. It is classically known that every Calder\'on-Zygmund operator $T$ extends to a bounded operator on $L^2(w)$ if and only if $w$ is an $A_2$  Muckenhoupt weight, namely, if and only if
\[
 [ w  ] _{A_2} \equiv \sup_{I} \left  \langle w \rangle_I
 \langle w^{-1}\right \rangle_I  < \infty,
\]
where the supremum is taken over all intervals $I$ and $\left \langle w \right \rangle_I \equiv \frac{1}{|I|} \int_I w(x) dx.$  In contrast, the question of the dependence of the operator norm of $T$ on $[w]_{A_2},$ called the \emph{$A_2$ conjecture}, remained open for decades. Mathematicians first bounded the norms of simpler operators including the Martingale transform, Hilbert transform, and more generally, dyadic shifts \cite{ lpr, pet07, wit00}.  A key strategy for proving the linear bound reduces general Calder\'on-Zygmund operators to simpler operators.   Using a refined method of decomposing Calder\'on-Zygmund operators as sums of dyadic shifts, Hyt\"onen resolved the $A_2$ conjecture in 2012 in \cite{h12}, showing
\[
\| T  \|_{L^2(w) \rightarrow L^2(w)} \lesssim [w]_{A_2}
\] 
for all Calder\'on-Zygmund operators $T$. More recently, Conde-Alonso and Rey, Lerner, Lacey \cite{CAR14, lac15, l13, l10} have developed simpler proofs of the $A_2$ conjecture by controlling Calder\'on-Zygmund operators using basic operators called \emph{sparse operators}, which easily satisfy the linear bound.

A basic tool when controlling these easier operators is the dyadic Carleson Embedding Theorem, which says:

\begin{thm}[Carleson Embedding Theorem] 
\label{thm:CET1} Let $\{a_I\}_{I \in \mathcal{D}}$ be a sequence of nonnegative numbers indexed by the grid of dyadic intervals $\mathcal{D}.$ Then 
\[ \sum_{I \in \mathcal{D}}  a_I \langle w^{\frac{1}{2}} f \rangle_I^2 \le C_1 \| f \|^2_{L^2} \quad\forall f\in L^2 \text{ if and only if } \frac{1}{|J|} \sum_{I: I \subseteq J} a_I \left \langle w \right \rangle_I^2 \le C_2 \left \langle w \right \rangle_J \quad \forall J \in \mathcal{D}. \]
Moreover, $C_2\leq C_1 \le 4 C_2.$
\end{thm}
We are interested in the development of these ideas in the matrix setting. 

\subsection{Matrix $A_2$ Conjecture.}   The relevant definitions are as follows: a $d\times d$ matrix-valued function $W(x)$ is a \emph{matrix weight} if its entries are locally integrable and if $W(x)$ is a positive definite matrix  for almost every $x \in \mathbb{R}$. Given a matrix weight $W$, one can define
\[
L^2(W) \equiv \left\{ f \in L^2(\mathbb{R}, \mathbb{C}^d): \| f\|^2_{L^2(W)} \equiv \int_{\mathbb{R}} \left \| W^{\frac{1}{2}}(x) f (x) \right \|^2 dx < \infty \right \}. 
\] 
Many scalar results have already been developed in this setting. For example, Nazarov-Treil \cite{nt97} and Volberg \cite{vol97} characterized the boundedness of Calder\'on-Zygmund operators on these matrix weighted $L^2$ spaces. Indeed, they showed that for each scalar Calder\'on-Zygmund operator $T$, the matrix Calder\'on-Zygmund operator $\tilde{T}$ defined on $L^2(\mathbb{R}, \mathbb{C}^d)$ by applying $T$ component-wise is bounded on $L^2(W)$ if $W$ is a matrix $A_2$ weight, namely if
\[
\big [ W \big ] _{A_2} \equiv \sup_{I} \left \| \langle W \rangle_I^{\frac{1}{2}}
 \langle W^{-1} \rangle_I^{\frac{1}{2}} \right \|^2 < \infty,
\]
where $\| \cdot \|$ denotes the norm of the matrix acting on $\mathbb{C}^d$. In this paper, we also use $\| \cdot \|$ to denote the norm of a vector in $\mathbb{C}^d$, but it will be clear from the context whether we are considering matrices or vectors.

In the interim, the study of operators on matrix-weighted spaces has received a great deal of attention. See \cite{bpw14, bw14, bow2001, cg01, gold03, IKP, lt07, nptv02, rou03, vol97}.  However, the question of the sharp dependence of the operators norms on $[W]_{A_2}$, termed the \emph{Matrix $A_2$ Conjecture} is still open. In \cite{bpw14}, the two authors with S. Petermichl showed that for the Hilbert transform $H$, 
\[ 
\| H  \|_{L^2(W) \rightarrow L^2(W)} \lesssim [W]_{A_2}^{\frac{3}{2}} \textnormal{log}\, [W]_{A_2},
\]
for all matrix $A_2$ weights $W$. Although this is the best known estimate, the bound is unlikely to be sharp. One key obstruction to obtaining better norm bounds is the following matrix version of Theorem \ref{thm:CET1}. This theorem is mentioned in \cite{IKP} and discussed in depth in \cite{bw14}:
\begin{thm} 
\label{thm:CET2} Let $W$ be an $A_2$ weight and let $\{A_I\}_{I\in\mathcal{D}}$ be a sequence of positive semi-definite $d \times d$ matrices. Then 
\[
\sum_{I\in\mathcal{D}} \left\langle A_I \left\langle W^{\frac{1}{2}}f\right\rangle_I, \left\langle W^{\frac{1}{2}} f\right\rangle_I\right\rangle_{\mathbb{C}^d} \le C_1 \left\Vert f\right\Vert_{L^2}^2 \  \text{ iff  } \ \
\frac{1}{|J|} \sum_{I: I \subseteq J} \left \langle W \right \rangle_I A_I  \left \langle W \right \rangle_I \le C_2  \left \langle W \right \rangle_J \ \ \ \forall J \in \mathcal{D},
\]
where $C_2\leq C_1 \lesssim C_2 [W]_{R_2} [W]_{A_2}.$
\end{thm}

Notice that this theorem holds only for $A_2$ weights and the guaranteed relationship between the testing constant $C_2$ and embedding constant $C_1$ is 
\[ C_1 \lesssim [W]_{A_2} [W]_{R_2} C_2. \]
Here, the notation $A \lesssim B$ indicates $A \le C B$, where $C$ is a constant typically depending on the dimension $d$. The term $[W]_{R_2}$ is a reverse H\"older constant, which is discussed in detail in \cite{bw14}. For our purposes, it suffices that $[W]_{R_2} \le c(d) [W]_{A_2}$, where $c(d)$ is a dimensional constant. The appearance of this (essentially) extra $[W]^2_{A_2}$ in Theorem \ref{thm:CET2} makes it very difficult to use this result to prove sharp bounds. 

The problem in removing this additional term is that no known scalar proofs of Theorem \ref{thm:CET1} appear to generalize to the matrix setting. Rather, Theorem \ref{thm:CET2} has only been proven via the following auxiliary Carleson Embedding Theorem, which was established by Treil-Volberg in \cite{vt97} and extended by Isralowitz-Pott-Kwon in \cite{IKP}.

\begin{thm} \label{thm:CET} Let $W$ be a matrix $A_2$ weight and let $\left\{A_I\right\}_{I\in\mathcal{D}}$ be a sequence of positive semi-definite $d\times d$ matrices.  Then 
\[
\sum_{I\in\mathcal{D}} \left\langle A_I \left\langle W^{\frac{1}{2}} f\right\rangle_I, \left\langle W^{\frac{1}{2}} f\right\rangle_I\right\rangle_{\mathbb{C}^d} \le C_1 \left\Vert f\right\Vert_{L^2}^2 \  \text{ if } \ 
\frac{1}{|J|} \sum_{I:I \subseteq J} \left \| \left \langle W \right \rangle^{\frac{1}{2}}_I A_I  \left \langle W \right \rangle^{\frac{1}{2}}_I \right \| \le C_2 \ \ \forall J \in \mathcal{D},
\]
where $C_1 \lesssim [W]_{R_2} C_2 \lesssim   [W]_{A_2} C_2$.
\end{thm}

\subsection{Summary of Results.}

In this paper, we offer two new proofs of Theorem \ref{thm:CET}. The only drawback is that these proofs give $C_1 \lesssim [W]_{A_2}^2 C_2$, which is not the optimal constant. Nevertheless, these new proofs show that this embedding theorem has close ties to both the boundedness of simple maximal functions and $H^1$-BMO duality in the matrix case. We end with an application of Theorem \ref{thm:CET} to the study of sparse operators, which shows that Theorem \ref{thm:CET} still allows one to prove new estimates for operators acting on matrix weighted $L^2$ spaces. Here is the overview of the paper.

In Section \ref{sec:maximal}, we fix a matrix weight $W$ and consider the maximal function $M_W$ defined by
\[ M_W f(x) \equiv \sup_{I \in \mathcal{D}: x \in I} \left \| \left \langle W \right \rangle_I^{-\frac{1}{2}} \left \langle W^{\frac{1}{2}} f \right \rangle_I \right \| \textbf{1}_I(x). \]
Using straightforward estimates, in Subsection \ref{subsec:maximal}, we bound this maximal function using an auxiliary maximal function studied by Christ and Goldberg in \cite{cg01}. Using a slight modification of the Christ-Goldberg arguments, we conclude in Corollary \ref{cor:Maximal} that 
\[ \| M_W \|_{L^2 \rightarrow L^2} \lesssim [W]_{A_2}, \]
for all matrix $A_2$ weights $W$. Then, in Subsection \ref{subsec:proof1}, we develop a new stopping-time argument that, when paired with the maximal function bound, yields an elegant proof of Theorem \ref{thm:CET}. 

In Section \ref{sec:sequences}, we consider the following two spaces of sequences of $d\times d$ matrices indexed by the dyadic intervals $\mathcal{D}:$ 
\[\begin{aligned}
\mathcal{S} &\equiv \left \{ \{S_I \}: S^2(x) \equiv \sum_{I \in \mathcal{D}} \| S_I \|^2 \frac{ 1_I(x)}{|I|} \in L^1(\mathbb{R}) \text{ and } \| \{S_I\}\|_{\mathcal{S}} \equiv \| S(x) \|_{L^1(\mathbb{R})}\right \}; \\
 \mathcal{T} &\equiv \left \{ \{ T_I \}: \| \{ T_I \}\|_{\mathcal{T}}^2 \equiv \sup_{J \in \mathcal{D}} \left \| \frac{1}{|J|} \sum_{I \subseteq J} T_I T_I^* \right \| < \infty \right \}.
 \end{aligned}
\] 
Notice that the sequences in $\mathcal{S}$ are related to functions whose square functions are in $L^1(\mathbb{R})$. This shows $\mathcal{S}$ is an $H^1$-type space and similarly, $\mathcal{T}$ is a space of $BMO$-type sequences. These spaces are similar to matrix analogues of scalar sequence spaces studied by Lee-Lin-Lin in \cite{lll09}. In Theorem \ref{thm:Duality}, we modify arguments of Lee-Lin-Lin to establish that each $\{T_I\} \in \mathcal{T}$ induces a linear functional on $\mathcal{S}.$ Then in Subsection \ref{subsec:proof2}, we use this duality relationship, paired with the maximal function bound, to provide another proof of Theorem \ref{thm:CET}.

Finally, in Section \ref{sec:sparse}, we consider an application of Theorem \ref{thm:CET}. Specifically, we say that an operator $S:L^2(\mathbb{R}, \mathbb{C}^d) \rightarrow L^2(\mathbb{R}, \mathbb{C}^d)$ is \emph{sparse} if 
\[ S f \equiv \sum_{I \in \mathfrak{S}} \left \langle f \right \rangle_I \textbf{1}_1(x), 
\]
where $\mathfrak{S} \subseteq \mathcal{D}$ is a collection of dyadic intervals satisfying the following sparsity condition: for each $I \in \mathfrak{S}$,
\begin{equation} \label{eqn:sparse} \sum_{J \in ch_\mathfrak{S}(I)} |J| \le \frac{1}{2} |I|, \end{equation}
where the sum is restricted to the $\mathfrak{S}$-children of $I$, namely the maximal elements of $\mathcal{S}$ that are strictly contained in $I.$   The constant $\frac{1}{2}$ is largely unimportant; any $0<c<1$ will work to define sparse families.  We use Theorem \ref{thm:CET} to establish the following result:

\begin{thm} \label{thm:sparse} Let $W$ be a $d \times d$ matrix $A_2$ weight and let $S$ be a sparse operator on $L^2(\mathbb{R}, \mathbb{C}^d).$ Then 
\[ \| S \|_{L^2(W) \rightarrow L^2(W)} \lesssim [W]_{A_2}^{\frac{3}{2}}. \]
\end{thm} 
Notice that this dependence of a sparse operator's norm on $[W]_{A_2}$ is better than the current known dependence of the Hilbert transform's norm on $[W]_{A_2}.$

Moreover, recall that in the scalar setting, sparse operators are used in \cite{CAR14,lac15,l13} to prove bounds for general Calder\'on-Zgymund operators. Unfortunately, those exact reduction arguments do not appear to generalize to the matrix setting. Nevertheless, it seems likely that alternate arguments will demonstrate close connections between general Calder\'on-Zgymund operators and sparse operators and at the least, allow one to obtain comparable norm bounds for a subclass of Calder\'on-Zgymund operators in the matrix setting.

\section{Theorem \ref{thm:CET} via Maximal Functions \& A Stopping-Time Argument} \label{sec:maximal}
\subsection{The Relevant Maximal Function} \label{subsec:maximal}

Recall the maximal function of interest: 
\[ M_W f(x) \equiv \sup_{I:x \in I} \left \| \left \langle W \right \rangle_I^{-\frac{1}{2}} \left \langle W^{\frac{1}{2}} f \right \rangle_I \right \| \textbf{1}_I(x), \]
where the supremum is taken over dyadic intervals. Although this is the maximal function most closely related to the matrix Carleson Embedding Theorem, it is difficult to bound directly. Instead, first notice that we have the following simple pointwise bound:
\[ 
\begin{aligned}
M_W f(x) &= \sup_{I:x \in I} \left \| \left \langle W \right \rangle_I^{-\frac{1}{2}}  \left \langle W^{-1} \right \rangle_I^{-\frac{1}{2}}  \left \langle W^{-1} \right \rangle_I^{\frac{1}{2}} \left \langle W^{\frac{1}{2}} f \right \rangle_I \right \|   \\
& \le \sup_{I:x \in I} \left \| \left \langle W \right \rangle_I^{-\frac{1}{2}}  \left \langle W^{-1} \right \rangle_I^{-\frac{1}{2}} \right \| \left \|  \left \langle W^{-1} \right \rangle_I^{\frac{1}{2}} \left \langle W^{\frac{1}{2}} f \right \rangle_I \right \| \\
&\le \sup_{I: x \in I} \left \|  \left \langle W^{-1} \right \rangle_I^{\frac{1}{2}} \left \langle W^{\frac{1}{2}} f \right \rangle_I \right \| \\
&= \sup_{I: x \in I} \left \|   \frac{1}{|I|} \int_I \left \langle W^{-1} \right \rangle_I^{\frac{1}{2}} W^{\frac{1}{2}}(y) f(y) \ dy  \right \| \\
&\le \sup_{I: x \in I} \frac{1}{|I|} \int_I \left \|  \left \langle W^{-1} \right \rangle_I^{\frac{1}{2}} W^{\frac{1}{2}}(y) f(y) \right \| \ dy.
\end{aligned}
\]
Here, we used the fact that for each $I\in \mathcal{D}$, it follows that $ \| \langle W \rangle_I^{-\frac{1}{2}}  \langle W^{-1}  \rangle_I^{-\frac{1}{2}}  \| \le 1.$ This inequality can be easily deduced from Corollary 3.3 in \cite{vt97}. 

This sequence of inequalities motivates the following definition: given a matrix weight $V,$ define the auxiliary maximal function $\widetilde{M}_V$ by
\[ \widetilde{M}_Vf(x) \equiv  \sup_{I: x \in I} \frac{1}{|I|} \int_I \left \|  \left \langle V \right \rangle_I^{\frac{1}{2}} V^{-\frac{1}{2}}(y) f(y) \right \| \ dy. \]
Then, as long as $W$ and $W^{-1}$ are both matrix weights, as is the case when $W$ is an $A_2$ matrix weight,  our previous arguments show that pointwise
\begin{equation} \label{eqn:maxbd} M_Wf(x) \le  \widetilde{M}_{W^{-1}}f(x). \end{equation}
 Christ-Goldberg actually studied this auxiliary maximal function $\widetilde{M}_V$ in \cite{cg01}. In their Lemma 2.2, they showed as long as $V$ is an $A_2$ matrix weight, then this maximal operator is bounded on $L^2(\mathbb{R}, \mathbb{C}^d).$ A close reading of their proof also reveals the following linear dependence of the operator norm on the $A_2$ characteristic of $V$.
 
\begin{thm} \label{thm:Maximal} If $V$ is a $d\times d$ matrix $A_2$ weight, then 
\[ \| \widetilde{M}_{V} \|_{L^2 \rightarrow L^2} \lesssim [V]_{A_2}, \]
where the implied constant depends on the dimension $d$.
\end{thm}

We include the following modified version of their arguments to track the exact dependence of the operator norm on $[V]_{A_2}:$ 

\begin{proof} Fix a matrix weight $V \in A_2$. We first establish the following inequality:
\[ \frac{1}{|I|} \int_I  \left \| V^{-\frac{1}{2}}(y) \left \langle V \right \rangle^{\frac{1}{2}}_I  \right \|^{2 + 2\epsilon} dy \lesssim [V]_{A_2}^{1+\epsilon}, \]
for some small $\epsilon >0$ and all $I \in \mathcal{D}.$ To obtain this, fix $e \in \mathbb{C}^d$ and  $I \in \mathcal{D}$. Recall that \cite[Lemma 3.5]{vt97} says that
\[  \left \langle V \right \rangle_I^{\frac{1}{2}} V^{-1}(x) \left \langle V \right \rangle_I^{\frac{1}{2}} \]
is also in $A_2$ with $A_2$ characteristic $[V]_{A_2}.$ Then, as in \cite[Lemma 1.5]{lt07}, it is not too difficult to show that the function
\[ \left \langle \left \langle V \right \rangle_I^{\frac{1}{2}} V^{-1}(x) \left \langle V \right \rangle_I^{\frac{1}{2}} e, e \right \rangle_{\mathbb{C}^d}\]
is a scalar $A_2$ weight with $A_2$ characteristic at most $[V]_{A_2}$ for each $e\in \mathbb{C}^d.$ Furthermore, as noted by Wittwer in \cite{wit00}, a careful reading of Coifman-Fefferman's proof of the reverse H\"older inequality for scalar $A_2$ weights in \cite{cf74}  shows that if $v$ is a scalar $A_2$ weight and $\epsilon \equiv \frac{c}{[v]_{A_2}}$ for a small-enough constant $c$, then
\[ \frac{1}{|I|} \int_I v^{1+\epsilon}(y) dy \le \left( \frac{2}{|I|} \int_I v(y) dy \right)^{1 + \epsilon}. \]
We will apply this to the scalar $A_2$ weights 
\[  \left \langle \left \langle V \right \rangle_I^{\frac{1}{2}} V^{-1}(x) \left \langle V \right \rangle_I^{\frac{1}{2}} e_i, e_i \right \rangle_{\mathbb{C}^d}, \]
where the $\{e_i\}_{i=1}^{d}$ are the standard unit normal vectors in $\mathbb{C}^d.$ Then, by equating norm and trace of positive definite matrices (up to a dimensional constant), one can compute \[ 
\begin{aligned}
\frac{1}{|I|} \int_I  \left \| V^{-\frac{1}{2}}(y) \left \langle V \right \rangle^{\frac{1}{2}}_I  \right \|^{2 + 2\epsilon} dy & = 
\frac{1}{|I|} \int_I  \left \|  \left \langle V \right \rangle_I^{\frac{1}{2}} V^{-1}(y) \left \langle V \right \rangle_I^{\frac{1}{2}} \right \|^{1 + \epsilon} dy \\
& \lesssim \frac{1}{|I|}   \int_I  \left( \text{Tr} \left (  \left \langle V \right \rangle_I^{\frac{1}{2}} V^{-1}(y) \left \langle V \right \rangle_I^{\frac{1}{2}} \right ) \right)^{1 + \epsilon} dy \\
& \lesssim \frac{1}{|I|}   \int_I  \displaystyle \text{max}_{1 \le i \le d}  \left \langle \left \langle V \right \rangle_I^{\frac{1}{2}} V^{-1}(y) \left \langle V \right \rangle_I^{\frac{1}{2}} e_i, e_i \right \rangle_{\mathbb{C}^d}^{1 + \epsilon} \ dy\\
 & \le  \sum_{i =1}^d  \frac{1}{|I|} \int_I \left \langle \left \langle V \right \rangle_I^{\frac{1}{2}} V^{-1}(y) \left \langle V \right \rangle_I^{\frac{1}{2}} e_i, e_i \right \rangle_{\mathbb{C}^d}^{1 + \epsilon} dy \\
&  \le   \sum_{i =1}^d \left( \frac{2}{|I|}    \int_I  \left \langle \left \langle V \right \rangle_I^{\frac{1}{2}} V^{-1}(y) \left \langle V \right \rangle_I^{\frac{1}{2}} e_i, e_i \right \rangle_{\mathbb{C}^d} dy \right)^{1 + \epsilon}\\
& \lesssim [V]_{A_2}^{1 +\epsilon}.
\end{aligned}
\]
Set $r = 2 + 2\epsilon$ and let $r' <2$ be conjugate to $r.$ Then, for $f \in L^2(\mathbb{R}, \mathbb{C}^d)$,
\[ 
\begin{aligned}
\widetilde{M}_{V} f(x) 
& = \sup_{I: x \in I} \frac{1}{|I|} \int_I \left \|  \left \langle V \right \rangle_I^{\frac{1}{2}} V^{-\frac{1}{2}}(y) f(y) \right \| \ dy \\
& \le  \left( \frac{1}{|I|} \int_I  \left \| V^{-\frac{1}{2}}(y) \left \langle V \right \rangle^{\frac{1}{2}}_I  \right \|^{r} dy \right)^ \frac{1}{r} \sup_{I: x \in I} \left( \frac{1}{|I|} \int_I  \| f(y)\|^{r'} dy \right)^{\frac{1}{r'}}\\
& \lesssim [V]_{A_2}  \left( M ( \|f \|^{r'})(x)  \right)^{\frac{1}{r'}}.
\end{aligned}
\]
Define $p = \frac{2}{r'}$. Then $p>1$ and so, the standard maximal function $M$ maps $L^p(\mathbb{R}, \mathbb{C}^d)$ to $ L^p(\mathbb{R}, \mathbb{C}^d)$. But, then 
\[ 
\begin{aligned}
\|   \widetilde{M}_{V} f \|_{L^2}^2 \lesssim [V]^2_{A_2}\left  \|  \left( M ( \|f \|^{r'})  \right)^{\frac{1}{r'}} \right \|^2_{L^2} =  [V]^2_{A_2}  \left \| M ( \|f \|^{r'}) \right \|_{L^p}^2 \lesssim C(d) [V]^2_{A_2} \|  f \|_{L^2}^2,
\end{aligned}
\]
using the maximal function bound. This immediately implies that
\[ \| \widetilde{M}_{V}  \|_{L^2 \rightarrow L^2} \lesssim  [V]_{A_2},\]
as desired. \end{proof}

The pointwise bound \eqref{eqn:maxbd} paired with Theorem \ref{thm:Maximal} immediately gives the following corollary:

\begin{cor} \label{cor:Maximal} If $V$ is a $d\times d$ matrix $A_2$ weight, then 
\[ \| M_{V} \|_{L^2 \rightarrow L^2} \lesssim [V]_{A_2}, \]
where the implied constant depends on the dimension $d$.
\end{cor} 

\subsection{Proof One of Theorem \ref{thm:CET}}  \label{subsec:proof1}

Our first proof of the matrix Carleson Embedding Theorem is motivated by a classical proof of the standard embedding result Theorem \ref{thm:CET1}, which uses a stopping time argument and maximal function bound.  Due to complications  involving matrix inequalities, there is no clear generalization of the scalar proof that gives a matrix version of the standard Carleson Embedding Theorem \ref{thm:CET1}.  However, there is a stopping-time proof of Theorem \ref{thm:CET} using the same types of arguments that appear in the scalar set-up. Here is the proof:

\begin{proof}
Fix $f \in L^2(\mathbb{R}, \mathbb{C}^d)$ and for each $k \in \mathbb{Z}$, let $\mathcal{J}_k$ be the set of maximal dyadic intervals $I$ with
\[ 2^{k-1} \le \left \| \left \langle W \right \rangle_I^{-\frac{1}{2}} \left \langle W^{\frac{1}{2}}f \right \rangle_I \right \| \le 2^{k}.\]
We say that $J \in \mathcal{J}_k^*$ if $k$ is the largest integer such that $J \subseteq I$ for some $I \in \mathcal{J}_k.$ Now recall that the related maximal function 
\[ M_{W} f(x) =  \sup_{I}  \left \|  \langle W \rangle_I^{-\frac{1}{2}}  \langle W^{\frac{1}{2}}f  \rangle_I  \right \| \textbf{1}_I(x)  \]
is in $L^2(\mathbb{R}, \mathbb{C}^d)$ by Corollary \ref{cor:Maximal}. Using this, it is clear that for each $J \in \mathcal{D}$, either $J \in \mathcal{J}_k^*$ for some $k \in \mathbb{Z}$ or $ \langle W^{\frac{1}{2}}f  \rangle_J $ is the zero vector. Now, as a way to further explore the relationship between this stopping-time set up and the maximal function, consider the function
\[ g(x) \equiv \sum_{k \in \mathbb{Z}} \sum_{I \in \mathcal{J}_k}  \left \| \left \langle W \right \rangle_I^{-\frac{1}{2}} \left \langle W^{\frac{1}{2}}f \right \rangle_I \right \| \textbf{1}_I(x).\]
Again, as $M_Wf \in L^2(\mathbb{R}, \mathbb{C}^d)$, we know that for almost every $x \in \mathbb{R}$, there is a largest $K \in \mathbb{Z}$ such that $x \in I \in \mathcal{J}_K.$ Then, we can conclude $M_W f(x) \approx 2^K$ and further, as each $\mathcal{J}_k$ is a disjoint collection of intervals,
\[ g(x) =  \sum_{k \in \mathbb{Z}: k \le K} \sum_{I \in \mathcal{J}_k}  \left \| \left \langle W \right \rangle_I^{-\frac{1}{2}} \left \langle W^{\frac{1}{2}}f \right \rangle_I \right \| \textbf{1}_I(x) 
\le  \sum_{k \in \mathbb{Z}: k \le K} 2^k \lesssim 2^{K} \lesssim M_W f(x). \]
Then this pointwise inequality and our definition of $\mathcal{J}_k$ gives
\[  \| M_W f \|_{L^2}^2 \gtrsim \| g \|_{L^2}^2 = \sum_{k \in \mathbb{Z}} \sum_{I \in \mathcal{J}_k}  \left \| \left \langle W \right \rangle_I^{-\frac{1}{2}} \left \langle W^{\frac{1}{2}}f \right \rangle_I \right \|^2 |I| \gtrsim \sum_{k \in \mathbb{Z}} 2^{2k}  \left | \bigcup_{I \in J_k} I \right |.\]
Now, assume the Carleson Embedding testing condition:
\[ \sum_{ J \subseteq I} \left \| \left \langle W \right \rangle^{\frac{1}{2}}_J A_J  \left \langle W \right \rangle^{\frac{1}{2}}_J  \right \| \le  C_2 |I|.\]
Then, we can compute
\[ 
\begin{aligned} \sum_{J \in \mathcal{D}} \left \langle A_J \left \langle W^{\frac{1}{2}} f \right \rangle_J,   \left \langle W^{\frac{1}{2}} f \right \rangle_J \right \rangle_{\mathbb{C}^d} & =  \sum_{k \in \mathbb{Z}}  \sum_{J \in \mathcal{J}_k^*} \left \langle A_J \left \langle W^{\frac{1}{2}} f \right \rangle_J,  \left \langle W^{\frac{1}{2}} f \right \rangle_J \right \rangle_{\mathbb{C}^d} \\
& \le \sum_{k \in \mathbb{Z}} \sum_{I \in \mathcal{J}_k} \sum_{\substack{J \subseteq I \\ J \in \mathcal{J}^*_k}}  \left \langle A_J \left \langle W^{\frac{1}{2}} f \right \rangle_J,   \left \langle W^{\frac{1}{2}} f \right \rangle_J \right \rangle_{\mathbb{C}^d} \\
& \le \sum_{k \in \mathbb{Z}}  \sum_{I \in \mathcal{J}_k}  \sum_{\substack{J \subseteq I \\ J \in \mathcal{J}^*_k}}
\left \| \left \langle W \right \rangle_J^{- \frac{1}{2}} \left \langle W^{\frac{1}{2}} f \right \rangle_J \right \|^2 \left \| \left \langle W \right \rangle_J^{\frac{1}{2}}  A_J \left \langle W \right \rangle_J^{ \frac{1}{2}}  \right \| \\
& \le  \sum_{k \in \mathbb{Z}}  2^{2k}  \sum_{I \in \mathcal{J}_k} \sum_{J \subseteq I} \left \| \left \langle W \right \rangle_J^{ \frac{1}{2}}  A_J \left \langle W \right \rangle_J^{\frac{1}{2}}  \right \| \\
& \le C_2 \sum_{k \in \mathbb{Z}}  2^{2k}  \sum_{I \in \mathcal{J}_k}   \left | I \right|\\
& = C_2 \sum_{k \in \mathbb{Z}}  2^{2k} \left | \bigcup_{I  \in \mathcal{J}_k} I \right| \\
&  \lesssim C_2 \left \| M_W f  \right \|^2_{L^2} \\
& \lesssim C_2 [W]^2_{A_2} \|f \|^2_{L^2},
\end{aligned}
\]
which gives the desired embedding result.
\end{proof}

\section{Theorem \ref{thm:CET} via $H^1$-BMO Duality and A Maximal Function} \label{sec:sequences}

In this section, we offer an alternate proof of Theorem \ref{thm:CET}. The idea is to prove the result via duality using a pairing between $H^1$-type sequences and $BMO$-type sequences. 

\subsection{Relevant Sequence Spaces} To establish the needed duality result, we study related spaces of sequences $\mathcal{S}$ and $\mathcal{T}$, which are composed of sequences of $d\times d$ matrices indexed by the dyadic intervals. As mentioned in the introduction, we study this $H^1$-type space of matrix sequences
\[
\mathcal{S} \equiv \left \{ \{S_I \}: S^2(x) \equiv \sum_{I \in \mathcal{D}} \| S_I \|^2 \frac{ 1_I(x)}{|I|} \in L^1(\mathbb{R}) \text{ and } \| \{S_I\}\|_{\mathcal{S}} \equiv \| S(x) \|_{L^1(\mathbb{R})}\right \}, \]
and this BMO-type space of matrix sequences
\[ \mathcal{T} \equiv \left \{ \{ T_I \}: \| \{ T_I \}\|^2_{\mathcal{T}} \equiv \sup_{J \in \mathcal{D}} \left \| \frac{1}{|J|} \sum_{I \subseteq J} T_I T_I^* \right \| < \infty \right \}.
\] 
Now, we modify the arguments of  Lee-Lin-Lin from \cite{lll09}, which were used to study different, but related scalar sequence spaces to obtain the following result. 
\begin{thm} \label{thm:Duality} For each $\{T_I\}$ in $\mathcal{T}$, the linear functional 
\[ \{S_I\} \mapsto \sum_{I\in \mathcal{D}} \text{Tr} \left( S_I T^*_I\right) \]
is continuous on $\mathcal{S}.$ Namely, there is a dimensional constant $c(d)$ (not depending on $\{ T_I \})$ such that for all $\{T_I\} \in \mathcal{T}$,
\begin{equation} \label{eqn:ss} \left | \sum_{I\in \mathcal{D}} \text{Tr} \left( S_I T^*_I\right) \right |  \le c(d) \| \{T_I\} \|_{\mathcal{T}} \|\{S_I\} \|_{\mathcal{S}} \qquad \forall \ \{S_I \} \in \mathcal{S}.\end{equation}
\end{thm} 

A complete analogue of the  Lee-Lin-Lin result from \cite{lll09} would also show that every continuous linear functional on $\mathcal{S}$ is induced by a sequence $\{T_I\}$ in $\mathcal{T}$. Although such a result is likely true in the context, we do not prove it because we do not require that to obtain Theorem \ref{thm:CET}.

\begin{proof} For a fixed $\{T_I\} \in \mathcal{T}$, we will establish \eqref{eqn:ss}. To this end, fix $\{S_I\} \in \mathcal{S}$ and for each $k \in \mathbb{Z}$, define the sets
\[ 
\begin{aligned}
\Omega_k &\equiv \left \{ x \in \mathbb{R}: S(x) > 2^k \right \}; \\
B_k & \equiv \left \{ I \in \mathcal{D}: |I \cap \Omega_k | > \frac{1}{2} | I | \ \text{ and } \ |I \cap \Omega_{k+1} |\le \frac{1}{2} |I| \right\},
\end{aligned}
\]
and let $\tilde{I}$ denote the maximal intervals in $B_k.$
Further, define the enlargement $\widetilde{\Omega}_k$ of $\Omega_k$ as follows:
\[ \widetilde{\Omega}_k \equiv \left \{ x \in \mathbb{R}: M( 1_{\Omega_k})(x) > \frac{1}{2} \right \}, \]
where $M$ denotes the Hardy-Littlewood maximal function. Furthermore, notice that
if $I \in B_k$ then $I \subset \widetilde{\Omega}_k.$ Further, as $S(x) \in L^1(\mathbb{R})$, observe that if $I \not \in B_k$ for every $k \in \mathbb{Z}$, then it must be the case that $S_I \equiv 0.$ Then, we can compute
\[
\begin{aligned}
\sum_{I\in \mathcal{D}} \text{Tr} \left( S_I T^*_I\right)  & =  \sum_{k \in \mathbb{Z}} \sum_{\tilde{I} \in B_k} \sum_{\substack{ I \subseteq \tilde{I} \\ I \in B_k}}  \text{Tr} \left( S_I T^*_I\right)  \\
&\le \sum_{k \in \mathbb{Z}} \sum_{\tilde{I} \in B_k}  \left( \sum_{\substack{ I \subseteq \tilde{I} \\ I \in B_k}} \|S_I \|^2 \right)^{\frac{1}{2}}  \left( \sum_{\substack{ I \subseteq \tilde{I} \\ I \in B_k}} \|T_I \|^2 \right)^{\frac{1}{2}} \\
& \le c(d) \sum_{k \in \mathbb{Z}} \sum_{\tilde{I} \in B_k}  \left( \sum_{\substack{ I \subseteq \tilde{I} \\ I \in B_k}} \|S_I \|^2 \right)^{\frac{1}{2}}  \left( \Bigg \|  \sum_{\substack{ I \subseteq \tilde{I}}} T_I T_I^*   \Bigg \|\right)^{\frac{1}{2}} \\
& \le c(d) \| \{T_I \} \|_{\mathcal{T}} \sum_{k \in \mathbb{Z}} \sum_{\tilde{I} \in B_k} | \tilde{I} |^{\frac{1}{2}}  
\left( \sum_{\substack{ I \subseteq \tilde{I} \\ I \in B_k}} \|S_I \|^2 \right)^{\frac{1}{2}}   \\
& \le c(d) \| \{T_I \} \|_{\mathcal{T}} \sum_{k \in \mathbb{Z}} |\widetilde{\Omega}_k |^{\frac{1}{2}}  \left( \sum_{ I \in B_k} \|S_I \|^2 \right)^{\frac{1}{2}},
\end{aligned}
\]
where we used Cauchy-Schwarz again and the fact that the $\tilde{I}$ in each $B_k$ are disjoint. Now we show:
\begin{equation} \label{eqn:Sest} \sum_{I \in B_k} \|S_I \|^2 \le 2^{2k+3} |\widetilde{\Omega}_k|. \end{equation}
First observe that
\[ \int_{\tilde{\Omega}_k \setminus \Omega_{k+1}}  S^2(x) \ dx \le 2^{2k+2} |\widetilde{\Omega}_k |, \]
using the definition of $\Omega_{k+1}$ and similarly
\[ 
\begin{aligned}
\int_{\tilde{\Omega}_k \setminus \Omega_{k+1}}  S^2(x)  &\ge \int _{\tilde{\Omega}_k \setminus \Omega_{k+1}}  \sum_{I \in B_K} \|S_I ||^2 \frac{1_I(x)}{|I|} \ dx \\
& =  \sum_{I \in B_k} \|S_I \|^2 \frac{ | I\cap (\widetilde{\Omega}_k \setminus \Omega_{k+1})|}{|I|} \\
& =  \sum_{I \in B_k} \|S_I \|^2 \frac{ | I \setminus \Omega_{k+1}|}{|I|} \\
& \ge \frac{1}{2} \sum_{I \in B_k} \|S_I \|^2. 
\end{aligned}
\]
Combining those two estimates gives \eqref{eqn:Sest}. Given  \eqref{eqn:Sest}, our previously-calculated inequality becomes:
\[
\begin{aligned}
 \left | \sum_{I\in \mathcal{D}} \text{Tr} \left( S_I T^*_I\right) \right | &\le c(d) \| \{T_I \} \|_{\mathcal{T}} \sum_{k \in \mathbb{Z}} |\widetilde{\Omega}_k |^{\frac{1}{2}}  \left( \sum_{ I \in B_k} \|S_I \|^2 \right)^{\frac{1}{2}}  \\
 &\le c(d) \| \{T_I \} \|_{\mathcal{T}} \sum_{k \in \mathbb{Z}} |\widetilde{\Omega}_k| 2^{k+2} \\
 & \lesssim c(d) \| \{T_I \} \|_{\mathcal{T}} \sum_{k \in \mathbb{Z}} |\Omega_k| 2^{k} \\
& \lesssim  c(d) \| \{T_I \} \|_{\mathcal{T}}  \| S(x) \|_{L^1(\mathbb{R})} \\
& = c(d)  \| \{T_I \} \|_{\mathcal{T}}  \| \{S_I \} \|_{\mathcal{S}}, 
\end{aligned}
\]
which is the desired inequality. \end{proof}

\subsection{Proof Two of Theorem \ref{thm:CET}} \label{subsec:proof2}

\begin{proof} Let $W$ be a $d\times d$ matrix $A_2$ weight and assume $\{A_I\}_{I \in \mathcal{D}}$ is a sequence of positive semidefinite $d\times d$ matrices satisfying the testing condition:
\[ \sum_{ I \subset J} \left \| \left \langle W \right \rangle^{\frac{1}{2}}_I A_I  \left \langle W \right \rangle^{\frac{1}{2}}_I  \right \| \le  C_2 |J|,\]
for all $J \in \mathcal{D}$. Then we can write:
\[
\begin{aligned}
 \sum_{I \in \mathcal{D}} \left \langle A_I \left \langle W^{\frac{1}{2}} f \right \rangle_I,   \left \langle W^{\frac{1}{2}} f \right \rangle_I \right \rangle_{\mathbb{C}^d}  
 &=  \sum_{I \in \mathcal{D}} \left \| A_I^{\frac{1}{2}}  \left \langle W^{\frac{1}{2}} f \right \rangle_I \right \|^2 \\
 & = \left \|  \left\{ A_I^{\frac{1}{2}}  \left \langle W^{\frac{1}{2}} f \right \rangle_I  \right \} \right \|^2_{\ell^2(\mathcal{D}, \mathbb{C}^d)}.
\end{aligned}
\]
We estimate this quantity via duality. Specifically, for each $\{ b_I \} \in \ell^2(\mathcal{D}, \mathbb{C}^d)$, we compute
\[ 
\begin{aligned}
\sum_{I \in \mathcal{D}} \left \langle A_I^{\frac{1}{2}}  \left \langle W^{\frac{1}{2}} f \right \rangle_I, b_I \right \rangle_{\mathbb{C}^d} 
&= \sum_{I \in \mathcal{D}} \text{Tr} \left(A_I^{\frac{1}{2}}  \left \langle W^{\frac{1}{2}} f \right \rangle_I b_I^* \right) \\
& = \sum_{I \in \mathcal{D}} \text{Tr} \left(A_I^{\frac{1}{2}}  \left \langle W \right \rangle_I^{\frac{1}{2}}   \left \langle W \right \rangle_I^{-\frac{1}{2}} \left \langle W^{\frac{1}{2}} f \right \rangle_I b_I^* \right). 
\end{aligned}
\]
To invoke Theorem \ref{thm:Duality}, define the sequences $\{S_I\}$ and $\{T_I\}$ by
\[ 
T_I =  \left \langle W \right \rangle_I^{\frac{1}{2}}A_I^{\frac{1}{2}} \text{ and } S_I = \left \langle W \right \rangle_I^{-\frac{1}{2}} \left \langle W^{\frac{1}{2}} f \right \rangle_I b_I^* \qquad \forall \ I \in \mathcal{D}.
\]
Notice that the testing condition implies that 
\[ \| \{ T_I \} \|^2_{\mathcal{T}} =  \sup_{J \in \mathcal{D}}  \left \| \frac{1}{|J|} \sum_{I \subseteq J} T_I T_I^* \right \| =  \sup_{J \in \mathcal{D}} \left \| \frac{1}{|J|} \sum_{I \subseteq J}  \left \langle W \right \rangle_I^{\frac{1}{2}}A_I  \left \langle W \right \rangle_I^{\frac{1}{2}}  \right  \| \le c(d) C_2. \]
Then Theorem \ref{thm:Duality} implies that
\[
\begin{aligned}
\sum_{I \in \mathcal{D}} \left \langle A_I^{\frac{1}{2}}  \left \langle W^{\frac{1}{2}} f \right \rangle_I, b_I \right \rangle_{\mathbb{C}^d} &= \sum_{I \in \mathcal{D}}  \text{Tr} \left( T^*_I S_I \right) \\
& = \sum_{I \in \mathcal{D}}  \text{Tr} \left( S_I T_I^* \right) \\ 
&\lesssim \left \| \{T_I \} \right \|_{\mathcal{T}}   \left \| \{ S_I \} \right \|_{\mathcal{S}} \\
& \lesssim \sqrt{C_2}  \left \| \{ S_I \}  \right\|_{\mathcal{S}}.
\end{aligned}
\]
Now observe that
\[ 
\begin{aligned}
S^2(x) &= \sum_{I \in \mathcal{D}} \|S_I\|^2 \frac{1_I(x)}{|I|} \\
& \le \sum_{I \in \mathcal{D}} \left \|  \left \langle W \right \rangle_I^{-\frac{1}{2}} \left \langle W^{\frac{1}{2}} f \right \rangle_I \right \|^2  \| b_I \|^2 \frac{1_I(x)}{|I|} \\
&\le \sup_{I: x \in I} \left \|  \left \langle W \right \rangle_I^{-\frac{1}{2}} \left \langle W^{\frac{1}{2}} f \right \rangle_I \right \|^2 \sum_{I \in \mathcal{D}}  \|b_I \|^2 \frac{1_I(x)}{|I|}.
\end{aligned}
\]
Notice that this is exactly the maximal function we studied earlier. Then Corollary \ref{cor:Maximal} implies that
\[ 
\begin{aligned}
\|S(x) \|_{L^1(\mathbb{R})}&\le \| M_W f \|_{L^2(\mathbb{R})} \left( \int_{\mathbb{R}} \sum_{I \in \mathcal{D}} \|b_I \|^2 \frac{1_I(x)}{|I|} dx \right)^{\frac{1}{2}} \\
&\lesssim [W]_{A_2} \|f \|_{L^2(\mathbb{R})} \| \{ b_I \} \|_{\ell^2(\mathcal{D}, \mathbb{C}^d)},
\end{aligned}
\]
as desired. So we can combine this with our previous estimates to conclude 
\[ \sum_{I \in \mathcal{D}} \left \langle A_I^{\frac{1}{2}}  \left \langle W^{\frac{1}{2}} f \right \rangle_I, b_I \right \rangle_{\mathbb{C}^d} \lesssim \sqrt{C_2}  [W]_{A_2} \|f \|_{L^2(\mathbb{R})} \| \{ b_I \} \|_{\ell^2(\mathcal{D}, \mathbb{C}^d)}. \] 
Since this estimate holds for all $\{b_I\} \in \ell^2(\mathcal{D}, \mathbb{C}^d)$, we immediately have
\[ \sum_{I \in \mathcal{D}}  \left \| A_I^{\frac{1}{2}}  \left \langle W^{\frac{1}{2}} f \right \rangle_I  \right \|^2 \lesssim C_2 [W]^2_{A_2} \| f \|^2_{L^2(\mathbb{R})}, \]
which completes the proof. \end{proof}

\section{Application: The Bound for Sparse Operators and the Proof of Theorem \ref{thm:sparse}} \label{sec:sparse}
Recall that an operator $S: L^2(\mathbb{R}, \mathbb{C}^d) \rightarrow L^2(\mathbb{R}, \mathbb{C}^d) $ is called \emph{sparse} if 
\[ S f = \sum_{I \in \mathfrak{S}} \left \langle f \right \rangle_I \textbf{1}_I, \]
where the collection of intervals $\mathfrak{S} \subseteq \mathcal{D}$ satisfies the sparseness condition given in \eqref{eqn:sparse}. Now we use Theorem \ref{thm:CET} to establish Theorem \ref{thm:sparse}, which basically says:
\[ 
\| S \|_{L^2(W) \rightarrow L^2(W)} \lesssim [W]_{A_2}^{\frac{3}{2}},\]
for every $d\times d$ matrix $A_2$ weight $W$. Here is the proof:

\begin{proof} Let $S$ be a sparse operator and observe that standard arguments give
\[ \| S \|_{L^2(W) \rightarrow L^2(W)}  = \| S M_{W^{-1}} \|_{L^2(W^{-1}) \rightarrow L^2(W)}.\]
So, we will study the second term instead and prove the desired bound using duality. Specifically, fix $f \in L^2(W^{-1})$ and $g \in L^2(W)$. Then 
\[
\begin{aligned}
\left \langle S M_{W^{-1}} f , g \right \rangle _{L^2(W)} & = \sum_{I \in \mathfrak{S}} \left \langle \left \langle W^{-1} f \right \rangle_I \textbf{1}_I, g \right \rangle_{L^2(W)} \\ 
&= \sum_{I \in \mathfrak{S}} \left \langle \left \langle W^{-1} f \right \rangle_I, \left \langle W g \right \rangle_I \right \rangle_{\mathbb{C}^d} |I| \\
&\le \sum_{I \in \mathfrak{S}} \left \|  \left \langle W  \right \rangle_I^{\frac{1}{2}} \left \langle W^{-1} \right \rangle_I^{\frac{1}{2}} \right \| \left \| \left \langle W^{-1}  \right \rangle_I^{-\frac{1}{2}}\left \langle W^{-1} f \right \rangle_I \right \| \left \|  \left \langle W \right \rangle_I^{-\frac{1}{2}} \left \langle W g \right \rangle_I \right \| | I| \\
&\le [W]_{A_2}^{\frac{1}{2}} \left( \sum_{I \in \mathfrak{S}}  \left \| \left \langle W^{-1}  \right \rangle_I^{-\frac{1}{2}}\left \langle W^{-1} f \right \rangle_I \right \|^2  |I| \right)^{\frac{1}{2}} \left( \sum_{I \in \mathfrak{S}} \left \|  \left \langle W \right \rangle_I^{-\frac{1}{2}} \left \langle W g \right \rangle_I \right \|^2 | I| \right)^{\frac{1}{2}}. 
\end {aligned} 
\]
We will show how to control the first sum above. The second will follow using symmetric arguments. First, observe that 
\[ \begin{aligned} 
\sum_{I \in \mathfrak{S}}  \left \| \left \langle W^{-1}  \right \rangle_I^{-\frac{1}{2}}\left \langle W^{-1} f \right \rangle_I \right \|^2  |I|  &= \sum_{I \in \mathfrak{S}} \left \langle  \left \langle W^{-1}  \right \rangle_I^{-1} |I| \left \langle W^{-1} f \right \rangle_I, \left \langle W^{-1} f \right \rangle_I\right \rangle_{\mathbb{C}^d} \\
&= \sum_{I \in \mathcal{D}} \left \langle A_I \left \langle W^{-1} f \right \rangle_I,  \left \langle W^{-1} f \right \rangle_I \right \rangle_{\mathbb{C}^d}, 
\end{aligned}
\]
where $\{A_I\}_{I \in \mathcal{D}}$ is the sequence of positive semidefinite matrices indexed by the dyadic intervals defined by:
\[ A_I = \left \{ \begin{array}{cc}
\left \langle W^{-1}  \right \rangle_I^{-1}|I| & \text{ if } I \in \mathfrak{S} \\
0 & \text{ if } I \not \in \mathfrak{S}. 
\end{array} \right.
\]
We wish to apply Theorem \ref{thm:CET} with the weight $W^{-1}.$ To do this, we must establish the appropriate testing conditions. Specifically, notice that if $J \in \mathcal{D}$, then
\[
\begin{aligned}
 \frac{1}{|J|} \sum_{I: I \subseteq J} \left \| \left \langle W^{-1} \right \rangle_I^{\frac{1}{2}}  A_I  \left \langle W^{-1} \right \rangle_I^{\frac{1}{2}}  \right \| &= \frac{ 1}{|J|} \sum_{I \subseteq J: I \in 
\mathfrak{S}}  \left \| \left \langle W^{-1} \right \rangle_I^{\frac{1}{2}}  \left \langle W^{-1}  \right \rangle_I^{-1}|I|  \left \langle W^{-1} \right \rangle_I^{\frac{1}{2}}  \right \|  \\
& = \frac{1}{|J|}  \sum_{I \subseteq J: I \in 
\mathfrak{S}} |I| \\
& \le \frac{1}{|J|}  \left( |J| + \frac{|J|}{2} + \frac{|J|}{4} + \dots \right)\\
& = 2,
\end{aligned}
\]
where we used the sparsity condition on $J$, the $\mathfrak{S}$-children of $J$, the $\mathfrak{S}$-children of the $\mathfrak{S}$-children of $J$, and so on. Now, by Theorem \ref{thm:CET}, we can conclude that
\[ \sum_{I \in \mathfrak{S}}  \left \| \left \langle W^{-1}  \right \rangle_I^{-\frac{1}{2}}\left \langle W^{-1} f \right \rangle_I \right \|^2  |I|  \lesssim [W]_{A_2} \| W^{-\frac{1}{2}} f \|^2_{L^2} = [W]_{A_2} \|f \|^2_{L^2(W^{-1})}. \]
We can similarly conclude that 
\[ \sum_{I \in \mathfrak{S}}  \left \| \left \langle W  \right \rangle_I^{-\frac{1}{2}}\left \langle W g \right \rangle_I \right \|^2  |I|  \lesssim [W]_{A_2} \| W^{\frac{1}{2}} g \|^2_{L^2} =  [W]_{A_2} \|g \|^2_{L^2(W)}. \]
It immediately follows that 
\[ \left \langle S M_{W^{-1}} f , g \right \rangle _{L^2(W)}  \lesssim [W]_{A_2}^{\frac{3}{2}} \| f \|_{L^2(W^{-1})} \| g \|_{L^2(W)}. \]
As $f$ and $g$ were arbitrary, we can conclude that 
\[ \| S \|_{L^2(W) \rightarrow L^2(W)}  = \| S M_{W^{-1}} \|_{L^2(W^{-1}) \rightarrow L^2(W)} \lesssim [W]_{A_2}^{\frac{3}{2}}, \]
as desired.\end{proof}


\end{document}